\documentclass [12pt]{article}
\usepackage{graphicx,amssymb,amsfonts,latexsym,amsmath,amsthm,times}
\usepackage{epsfig}
\usepackage{color}
\usepackage[english]{babel}
\usepackage[a4paper]{geometry}
\def\lanbox{\hbox{$\, \vrule height 0.25cm width 0.25cm depth 0.01cm \,$}}

\numberwithin{equation}{section}

\begin{document}

\vspace*{1.4cm}

\normalsize \centerline{\Large \bf ON THE WELL-POSEDNESS
OF A CERTAIN MODEL WITH}

\medskip

\centerline{\Large \bf TWO KERNELS APPEARING IN THE MATHEMATICAL BIOLOGY}

\vspace*{1cm}

\centerline{\bf Messoud Efendiev$^{1,2}$, Vitali Vougalter$^{3 \ *}$}

\bigskip

\centerline{$^1$ Helmholtz Zentrum M\"unchen, Institut f\"ur Computational
Biology, Ingolst\"adter Landstrasse 1}

\centerline{Neuherberg, 85764, Germany}

\centerline{e-mail: messoud.efendiyev@helmholtz-muenchen.de}

\centerline{$^2$ Azerbaijan University of Architecture and Construction, Baku, Azerbaijan}

\centerline{e-mail: messoud.efendiyev@gmail.com}

\bigskip

\centerline{$^{3 \ *}$  Department of Mathematics, University
of Toronto}

\centerline{Toronto, Ontario, M5S 2E4, Canada}

\centerline{ e-mail: vitali@math.toronto.edu}

\medskip

%*******************************************************************
%ABSTRACT
%*******************************************************************

\vspace*{0.25cm}

\noindent {\bf Abstract:}
The work is devoted to establishing the global well-posedness in
$W^{(1,2),2}({\mathbb R}\times {\mathbb R}^{+})$ of the
integro-differential problem involving the two nonlocal terms describing the 
diffusion and the production in the biological system in
the presence of the transport term. Such model is relevant to the cell
population dynamics in the Mathematical Biology.
The proof is based on a fixed point technique.

\vspace*{0.25cm}



\vspace*{0.5cm}

\bigskip

\bigskip

%%%%%%%%%%%%%%%%%%%%%%%%%%%%%%%%%%%%%%%%%%%%%%%%%%

\setcounter{section}{1}

\centerline{\bf 1. Introduction}

\medskip

In the present article we establish  the global well-posedness
of the nonlocal reaction-diffusion equation with the constants
$a\geq 0$ and $b\in {\mathbb R}$, namely
\begin{equation}
\label{h}
\frac{\partial u}{\partial t} =
\int_{-\infty}^{\infty}J(x-y)u(y,t)dy +
b\frac{\partial u}{\partial x}+au+
\int_{-\infty}^{\infty}G(x-y)F(u(y,t), y)dy, 
\end{equation}
relevant to the cell population dynamics. We assume that the initial
condition for (\ref{h}) is 
\begin{equation}
\label{ic}
u(x,0)=u_{0}(x)\in H^{2}({\mathbb R}).
\end{equation}
The  global well-posedness of the integro-differential problem analogous to 
(\ref{h}) involving the fractional Laplacian in the diffusion term in the
context of the anomalous diffusion was discussed in ~\cite{EV25}.
The case of the bi-Laplacian in the diffusion term was covered in ~\cite{EV251}.
Spatial structures and generalized travelling waves for an
integro-differential equation were treated in ~\cite{ABVV10}. Article
~\cite{VV130} is devoted to the 
emergence and propagation of patterns in nonlocal reaction-
diffusion equations arising in the theory of speciation and containing the
transport term.  In ~\cite{GVA06} the authors consider the pattern and waves for a model in population 
dynamics with nonlocal consumptions of resources.
The existence of steady states and travelling waves for the
non-local Fisher-KPP equation was established in ~\cite{BNPR09}. In
~\cite{BHN05} the authors estimated the speed of propagation for KPP type
problems in the periodic framework. Important applications to the theory of
reaction-diffusion equations with non-Fredholm operators were developed in
~\cite{DMV05}, ~\cite{DMV08}. Entropy method for generalized Poisson-Nernst-
Planck equations was developed in ~\cite{GK18}. Work ~\cite{GHLP22} deals with
the local and global existence for nonlocal multispecies advection-diffusion models.
Persistence and time periodic positive solutions of doubly nonlocal Fisher-KPP equations in time periodic and space heterogeneous 
media were studied in ~\cite{GGS21}. Doubly nonlocal reaction-diffusion equations and the emergence of species were discussed in
~\cite{BVV17}. Lower estimate of the attractor dimension for a chemotaxis growth system was derived in ~\cite{ATEYM06}.
Quasilinear elliptic equations on half- and quarter-spaces were considered in ~\cite{DDE13}.
Exponential decay toward equilibrium via entropy methods for reaction-diffusion problems was established in ~\cite{DF06}.
Fredholm structures, topological invariants and applications were treated in ~\cite{E09}. Article ~\cite{E10} deals with
the development of the theory of finite and infinite dimensional attractors
for evolution equations of mathematical physics. Evolution equations arising in the modelling of life sciences were studied in 
~\cite{E13}. Work ~\cite{E131} is devoted to the attractors for degenerate parabolic type problems.
Mathematical analysis of a PDE--ODE coupled model of mitochondrial swelling with degenerate calcium ion diffusion was
performed in ~\cite{EOE20}.  The large time behavior of solutions of fourth order parabolic equations and $\epsilon$-entropy of their attractors were investigated in ~\cite{EP07}.  Spatial patterns arising in higher order models in physics and mechanics were covered
in ~\cite{PT01}.

The space variable $x$ in our work corresponds to the cell genotype,
$u(x,t)$ denotes the cell density as a function of the genotype and time.
The right side of (\ref{h}) describes the evolution of the cell density by virtue of
the cell proliferation, mutations and transport.
The nonlocal diffusion term here is correspondent to the change of the genotype
due to the small random mutations, and the integral nonlinear production term describes large mutations.
The function $F(u, x)$ stands for the rate of the cell birth dependent on $u$ and $x$
(density dependent proliferation), and the kernel $G(x-y)$ gives
the proportion of newly born cells changing their genotype from $y$ to $x$.
We assume that it depends on the distance between the genotypes.

The standard Fourier transform used throughout the present work is given by
\begin{equation}
\label{ft}
\widehat{\phi}(p)=\frac{1}{\sqrt{2\pi}}\int_{-\infty}^{\infty}\phi(x)e^{-ipx}dx,
\quad p\in {\mathbb R}.
\end{equation}
Clearly, the inequality
\begin{equation}
\label{fub}  
\|\widehat{\phi}(p)\|_{L^{\infty}({\mathbb R})}\leq \frac{1}{\sqrt{2\pi}}
\|\phi(x)\|_{L^{1}({\mathbb R})}
\end{equation}
is valid (see e.g. ~\cite{LL97}). Evidently, (\ref{fub}) implies
\begin{equation}
\label{fub1}  
\|p^{2}\widehat{\phi}(p)\|_{L^{\infty}({\mathbb R})}\leq \frac{1}{\sqrt{2\pi}}
\Big\|\frac{d^{2}\phi}{dx^{2}}\Big\|_{L^{1}({\mathbb R})}.
\end{equation}
We suppose that the conditions below on the integral kernels involved in 
problem (\ref{h}) are fulfilled.

\medskip

\noindent
{\bf Assumption 1.1.}  {\it Let $J(x), \ G(x): {\mathbb R}\to {\mathbb R}$ be
nontrivial, such that
$$
J(x), \ G(x), \ \frac{d^{2}G(x)}{dx^{2}}\in L^{1}({\mathbb R}).
$$
Moreover, the real part }
\begin{equation}
\label{re}
Re [\widehat{J}(p)]\leq 0, \ p\in{\mathbb R}. 
\end{equation}

\medskip

This allows us to define the auxiliary quantity
\begin{equation}
\label{q}  
\nu:=\sqrt{\|G(x)\|_{L^{1}({\mathbb R})}^{2}+
\Big\|\frac{d^{2}G(x)}{dx^{2}}\Big\|_{L^{1}({\mathbb R})}^{2}}.
\end{equation}
Thus, $0<\nu<\infty$.
  
\medskip

From the point of view of the applications, the space dimension is
not restricted to $d=1$ since the space variable corresponds to the cell
genotype but not to the usual physical space.
We have the Sobolev space
\begin{equation}
\label{ss}  
H^{2}({\mathbb R}):=
\Bigg\{\phi(x):{\mathbb R}\to {\mathbb {\mathbb R}} \ | \
\phi(x)\in L^{2}({\mathbb R}), \ \frac{d^{2}\phi}{dx^{2}}\in
L^{2}({\mathbb R}) \Bigg \}. 
\end{equation}
It is equipped with the norm
\begin{equation}
\label{n1}
\|\phi\|_{H^{2}({\mathbb R})}^{2}:=\|\phi\|_{L^{2}({\mathbb R})}^{2}+
\Bigg\| \frac{d^{2}\phi}{dx^{2}}\Bigg\|_{L^{2}({\mathbb R})}^{2}.
\end{equation}
Let us establish the global well-posedness for our problem (\ref{h}), (\ref{ic}).
For that purpose, we will use the function space
$$
W^{(1, 2), 2}({\mathbb R}\times [0, T]):=
$$
\begin{equation}
\label{122}
\Big\{u(x,t): {\mathbb R}\times [0, T]\to {\mathbb R} \ \Big|
\ u(x,t), \ \frac{\partial^{2}u}{\partial x^{2}}, \
\frac{\partial u}{\partial t}
\in L^{2}({\mathbb R}\times [0, T]) \Big\},
\end{equation}  
so that
$$
\|u(x,t)\|_{W^{(1, 2), 2}({\mathbb R}\times [0, T])}^{2}:=
$$
\begin{equation}
\label{122n}
\Big\|\frac{\partial u}{\partial t}\Big\|_{L^{2}({\mathbb R}\times [0, T])}^{2}+
\Big\|\frac{\partial^{2}u}{\partial x^{2}}\Big\|_{L^{2}({\mathbb R}\times [0, T])}^{2}+
\|u\|_{L^{2}({\mathbb R}\times [0, T])}^{2},
\end{equation}  
where $T>0$. Here
$$
\|u\|_{L^{2}({\mathbb R}\times [0, T])}^{2}:=\int_{0}^{T}\int_{-\infty}^{\infty}
|u(x,t)|^{2}dxdt.
$$
Throughout the work we will also use another norm
$$
\|u(x,t)\|_{L^{2}({\mathbb R})}^{2}:=\int_{-\infty}^{\infty}|u(x,t)|^{2}dx.
$$

\medskip

\noindent
{\bf Assumption 1.2.} {\it Function
$F(u, x): {\mathbb R}\times{\mathbb R}\to {\mathbb R}$ is satisfying the
Caratheodory condition (see ~\cite{K64}), such that
\begin{equation}
\label{lub}
|F(u, x)|\leq k|u|+h(x) \quad for \quad u\in {\mathbb R}, \quad x\in {\mathbb R}
\end{equation}  
with a constant $k>0$ and
$h(x): {\mathbb R}\to {\mathbb R}^{+}, \ h(x)\in L^{2}({\mathbb R})$.
Moreover, it is a Lipschitz continuous function, so that
\begin{equation}
\label{lip}
|F(u_{1}, x)-F(u_{2}, x)|\leq l|u_{1}-u_{2}| \quad for \quad any \quad
u_{1, 2}\in {\mathbb R}, \quad x\in {\mathbb R}
\end{equation}
with a constant $l>0$.}

\medskip

In the article ${\mathbb R}^{+}$ designates the nonnegative semi-axis.
The solvability of a local elliptic equation in a bounded domain in
${\mathbb R}^{N}$ was covered in ~\cite{BO86}. The nonlinear function involved there
was allowed to have a sublinear growth.
We apply the standard Fourier transform (\ref{ft}) to both sides of
problem (\ref{h}), (\ref{ic}) and arrive at
\begin{equation}
\label{hf}
\frac{\partial \widehat{u}}{\partial t}=[\sqrt{2\pi}\widehat{J}(p)+ibp+a]\widehat{u}+
\sqrt{2\pi}\widehat{G}(p)\widehat{f}_{u}(p,t),
\end{equation}
\begin{equation}
\label{icf}
\widehat{u}(p,0)=\widehat{u_{0}}(p).  
\end{equation}
Here and below $\widehat{f}_{u}(p,t)$ will denote the Fourier image
of $F(u(x,t), x)$. Obviously, we have
$$
u(x,t)=\frac{1}{\sqrt{2\pi}}\int_{-\infty}^{\infty}\widehat{u}(p,t)e^{ipx}dp,
\quad 
\frac{\partial u}{\partial t}=\frac{1}{\sqrt{2\pi}}\int_{-\infty}^{\infty}
\frac{\partial \widehat{u}(p,t)}{\partial t}e^{ipx}dp,
$$
where $x\in {\mathbb R}, \ t\geq 0$. By virtue of the Duhamel's principle,
we can reformulate problem (\ref{hf}), (\ref{icf}) as
$$
\widehat{u}(p,t)=
$$
\begin{equation}
\label{duh}
e^{t\{\sqrt{2\pi}\widehat{J}(p)+ibp+a\}}\widehat{u_{0}}(p)+\int_{0}^{t}
e^{(t-s)\{\sqrt{2\pi}\widehat{J}(p)+ibp+a\}}\sqrt{2\pi}\widehat{G}(p)\widehat{f}_{u}(p,s)ds.
\end{equation}
Let us consider the auxiliary problem related to equation (\ref{duh}), namely
$$
\widehat{u}(p,t)=
$$
\begin{equation}
\label{aux}
e^{t\{\sqrt{2\pi}\widehat{J}(p)+ibp+a\}}\widehat{u_{0}}(p)+\int_{0}^{t}
e^{(t-s)\{\sqrt{2\pi}\widehat{J}(p)+ibp+a\}}\sqrt{2\pi}\widehat{G}(p)\widehat{f}_{v}(p,s)ds.
\end{equation}
Here $\widehat{f}_{v}(p,s)$ stands for the Fourier image of
$F(v(x,s), x)$ under transform (\ref{ft}) and $v(x,t)\in W^{(1, 2), 2}({\mathbb R}\times [0, T])$ is arbitrary.

We introduce the operator $\tau_{a, b}$, so that $u =\tau_{a, b}v$, where $u$ 
satisfies (\ref{aux}). The main statement of the article is as follows.

\bigskip

\noindent
{\bf Theorem 1.3.} {\it Let Assumptions 1.1 and 1.2 be valid and
\begin{equation}
\label{qlt}
\nu l\sqrt{T^{2}e^{2aT}(1+2[a+|b|+\|J(x)\| _{L^{1}({\mathbb R})}]^{2})+2}<1
\end{equation}
with the constant $\nu$ introduced in (\ref{q}) and the Lipschitz constant
$l$ defined in (\ref{lip}).
Then equation (\ref{aux}) defines the map
$\tau_{a, b}: W^{(1, 2), 2}({\mathbb R}\times [0, T])\to W^{(1, 2), 2}
({\mathbb R}\times [0, T])$, which is a strict contraction.
The unique fixed point $w(x,t)$ of this map $\tau_{a, b}$ is the only solution of
problem (\ref{h}), (\ref{ic}) in
$W^{(1, 2), 2}({\mathbb R}\times [0, T])$.}

\medskip

The final proposition of the work is devoted to the global well-posedness for
our equation. 

\medskip

\noindent
{\bf Corollary 1.4.} {\it Let the conditions of Theorem 1.3 hold.
Then problem (\ref{h}), (\ref{ic}) has a unique solution   
$w(x,t)\in W^{(1, 2), 2}({\mathbb R}\times {\mathbb R}^{+})$. This solution 
does not vanish identically for $x\in {\mathbb R}$ and $t\in {\mathbb R}^{+}$ 
provided the intersection of supports of the Fourier images of functions
$\hbox{supp}\widehat{F(0, x)}\cap \hbox{supp}\widehat{G}$ is a set of
nonzero Lebesgue measure on the real line.}

\medskip

Let us turn our attention to the proof of our main statement.

\bigskip

%%%%%%%%%%%%%%%%%%%%%%%%%%%%%%%%%%%%%%%%%%%%%%%%%%

\setcounter{section}{2}
\setcounter{equation}{0}

\centerline{\bf 2. The well-posedness of the model}

\bigskip

\noindent
{\it Proof of Theorem 1.3.} We choose arbitrarily
$v(x,t)\in W^{(1, 2), 2}({\mathbb R}\times [0, T])$ and recall condition (\ref{re}).
Let us demonstrate that the first term in the right side of (\ref{aux})
is contained in $L^{2}({\mathbb R}\times [0, T])$. Clearly,
$$
\|e^{t\{\sqrt{2\pi}\widehat{J}(p)+ibp+a\}}\widehat{u_{0}}(p)\|_{L^{2}({\mathbb R})}^{2}=
\int_{-\infty}^{\infty}e^{2t\sqrt{2\pi}Re[\widehat{J}(p)]}e^{2at}|\widehat{u_{0}}(p)|^{2}dp\leq
e^{2at}\|u_{0}\|_{L^{2}({\mathbb R})}^{2},
$$
so that
$$
\|e^{t\{\sqrt{2\pi}\widehat{J}(p)+ibp+a\}}\widehat{u_{0}}(p)\|_{L^{2}({\mathbb R}\times [0, T])}^{2}=
\int_{0}^{T}\|e^{t\{\sqrt{2\pi}\widehat{J}(p)+ibp+a\}}\widehat{u_{0}}(p)\|_{L^{2}({\mathbb R})}^{2}dt
\leq
$$
$$
\int_{0}^{T}e^{2at}\|u_{0}\|_{L^{2}({\mathbb R})}^{2}dt.
$$
Obviously, this is equals to
$\displaystyle{\frac{e^{2aT}-1}{2a}\|u_{0}\|_{L^{2}({\mathbb R})}^{2}}$
when $a>0$ and
$T\|u_{0}\|_{L^{2}({\mathbb R})}^{2}$ for $a=0$. Thus,
\begin{equation}
\label{u0l2}  
e^{t\{\sqrt{2\pi}\widehat{J}(p)+ibp+a\}}\widehat{u_{0}}(p)\in L^{2}({\mathbb R}\times [0, T]).
\end{equation}
We derive the upper bound on the norm of the second term in the right side of (\ref{aux}) as
$$
\Big\|\int_{0}^{t}
e^{(t-s)\{\sqrt{2\pi}\widehat{J}(p)+ibp+a\}}\sqrt{2\pi}\widehat{G}(p)\widehat{f}_{v}(p,s)ds\Big\|_
{L^{2}({\mathbb R})}\leq
$$
$$  
\int_{0}^{t}\Big\|
e^{(t-s)\{\sqrt{2\pi}\widehat{J}(p)+ibp+a\}}\sqrt{2\pi}\widehat{G}(p)\widehat{f}_{v}(p,s)\Big\|_
{L^{2}({\mathbb R})}ds.
$$
Evidently,
$$  
\Big\|
e^{(t-s)\{\sqrt{2\pi}\widehat{J}(p)+ibp+a\}}\sqrt{2\pi}\widehat{G}(p)\widehat{f}_{v}(p,s)\Big\|_
{L^{2}({\mathbb R})}^{2}=
$$
\begin{equation}
\label{intgfv}  
\int_{-\infty}^{\infty}e^{2(t-s)\sqrt{2\pi}Re[\widehat{J}(p)]}e^{2a(t-s)}2\pi|\widehat{G}(p)|^{2}
|\widehat{f}_{v}(p,s)|^{2}dp.
\end{equation}
Let us use inequality (\ref{fub}) to obtain the estimate from above on the right side of
(\ref{intgfv}) as
$$
e^{2a(t-s)}2\pi\|\widehat{G}(p)\|_{L^{\infty}({\mathbb R})}^{2}
\|F(v(x,s),x)\|_{L^{2}({\mathbb R})}^{2}\leq
$$
$$
e^{2aT}\|G(x)\|_{L^{1}({\mathbb R})}^{2}
\|F(v(x,s),x)\|_{L^{2}({\mathbb R})}^{2}.
$$
Hence,
$$  
\Big\|
e^{(t-s)\{\sqrt{2\pi}\widehat{J}(p)+ibp+a\}}\sqrt{2\pi}\widehat{G}(p)\widehat{f}_{v}(p,s)\Big\|_
{L^{2}({\mathbb R})}\leq
$$
$$
e^{aT}\|G(x)\|_{L^{1}({\mathbb R})}
\|F(v(x,s),x)\|_{L^{2}({\mathbb R})}.
$$
By virtue of condition (\ref{lub}), we have
\begin{equation}
\label{lubs}  
\|F(v(x,s),x)\|_{L^{2}({\mathbb R})}\leq k\|v(x,s)\|_{L^{2}({\mathbb R})}+\|h(x)\|_
{L^{2}({\mathbb R})}.
\end{equation}
This means that 
$$  
\Big\|
e^{(t-s)\{\sqrt{2\pi}\widehat{J}(p)+ibp+a\}}\sqrt{2\pi}\widehat{G}(p)\widehat{f}_{v}(p,s)\Big\|_
{L^{2}({\mathbb R})}\leq
$$
$$
e^{aT}\|G(x)\|_{L^{1}({\mathbb R})}
\{k\|v(x,s)\|_{L^{2}({\mathbb R})}+\|h(x)\|_{L^{2}({\mathbb R})}\},
$$
so that
$$  
\Big\|\int_{0}^{t}
e^{(t-s)\{\sqrt{2\pi}\widehat{J}(p)+ibp+a\}}\sqrt{2\pi}\widehat{G}(p)\widehat{f}_{v}(p,s)ds\Big\|_
{L^{2}({\mathbb R})}\leq
$$
$$
ke^{aT}\|G(x)\|_{L^{1}({\mathbb R})}\int_{0}^{T}\|v(x,s)\|_{L^{2}({\mathbb R})}ds+
Te^{aT}\|G(x)\|_{L^{1}({\mathbb R})}\|h(x)\|_{L^{2}({\mathbb R})}.
$$
Let us use the Schwarz inequality, namely
\begin{equation}
\label{sch}  
\int_{0}^{T}\|v(x,s)\|_{L^{2}({\mathbb R})}ds\leq
\sqrt{\int_{0}^{T}\|v(x,s)\|_{L^{2}({\mathbb R})}^{2}ds}\sqrt{T}.
\end{equation}
This gives us
$$
\Big\|\int_{0}^{t}
e^{(t-s)\{\sqrt{2\pi}\widehat{J}(p)+ibp+a\}}\sqrt{2\pi}\widehat{G}(p)\widehat{f}_{v}(p,s)ds\Big\|_
{L^{2}({\mathbb R})}^{2}\leq
$$
$$
e^{2aT}\|G(x)\|_{L^{1}({\mathbb R})}^{2}
\{k\sqrt{T}\|v(x,s)\|_{L^{2}({\mathbb R}\times [0,T])}+T\|h(x)\|_{L^{2}({\mathbb R})}\}^{2}.
$$
We obtain the upper bound on the norm as
$$
\Big\|\int_{0}^{t}
e^{(t-s)\{\sqrt{2\pi}\widehat{J}(p)+ibp+a\}}\sqrt{2\pi}\widehat{G}(p)\widehat{f}_{v}(p,s)ds\Big\|_
{L^{2}({\mathbb R}\times [0,T])}^{2}=
$$
$$
\int_{0}^{T}\Big\|\int_{0}^{t}
e^{(t-s)\{\sqrt{2\pi}\widehat{J}(p)+ibp+a\}}\sqrt{2\pi}\widehat{G}(p)\widehat{f}_{v}(p,s)ds\Big\|_
{L^{2}({\mathbb R})}^{2}dt\leq 
$$
$$
e^{2aT}\|G(x)\|_{L^{1}({\mathbb R})}^{2}
\{k\|v(x,s)\|_{L^{2}({\mathbb R}\times [0,T])}+\sqrt{T}\|h(x)\|_{L^{2}({\mathbb R})}\}^{2}
T^{2}<\infty
$$
under the stated assumptions for
$v(x,s)\in W^{(1, 2), 2}({\mathbb R}\times [0, T])$. Thus,
\begin{equation}
\label{0tetsgvl2}  
\int_{0}^{t}
e^{(t-s)\{\sqrt{2\pi}\widehat{J}(p)+ibp+a\}}\sqrt{2\pi}\widehat{G}(p)\widehat{f}_{v}(p,s)ds\in
L^{2}({\mathbb R}\times [0,T]).
\end{equation}
Let us recall (\ref{aux}) along with statements (\ref{u0l2}) and  (\ref{0tetsgvl2}). Therefore,
\begin{equation}
\label{upt12}  
\widehat{u}(p,t)\in L^{2}({\mathbb R}\times [0,T]).
\end{equation}
Clearly,
\begin{equation}
\label{uxtl2}  
u(x,t)\in L^{2}({\mathbb R}\times [0,T]).
\end{equation}
By means of  (\ref{aux}), we have
$$
p^{2}\widehat{u}(p,t)=
$$
\begin{equation}
\label{aux2}
e^{t\{\sqrt{2\pi}\widehat{J}(p)+ibp+a\}}p^{2}\widehat{u_{0}}(p)+\int_{0}^{t}
e^{(t-s)\{\sqrt{2\pi}\widehat{J}(p)+ibp+a\}}\sqrt{2\pi}p^{2}\widehat{G}(p)\widehat{f}_{v}(p,s)ds.
\end{equation}
We use (\ref{re}) to treat the first term in the right side of (\ref{aux2}), so that
$$
\|e^{t\{\sqrt{2\pi}\widehat{J}(p)+ibp+a\}}p^{2}\widehat{u_{0}}(p)\|_{L^{2}({\mathbb R}\times [0,T])}^{2}=
\int_{0}^{T}\int_{-\infty}^{\infty}e^{2t\sqrt{2\pi}Re[\widehat{J}(p)]}e^{2at}|p^{2}\widehat{u_{0}}(p)|^{2}dp
dt\leq
$$
$$
\int_{0}^{T}\int_{-\infty}^{\infty}e^{2at}|p^{2}\widehat{u_{0}}(p)|^{2}dpdt.
$$
Obviously, this equals to
$\displaystyle
{\frac{e^{2aT}-1}{2a}\Big\|\frac{d^{2}u_{0}}{dx^{2}}\Big\|_{L^{2}({\mathbb R})}^{2}}$ for
$a>0$ and
$\displaystyle{T\Big\|\frac{d^{2}u_{0}}{dx^{2}}\Big\|_{L^{2}({\mathbb R})}^{2}}$
if $a=0$. This yields
\begin{equation}
\label{p2u0hpl2}
e^{t\{\sqrt{2\pi}\widehat{J}(p)+ibp+a\}}p^{2}\widehat{u_{0}}(p)\in L^{2}({\mathbb R}\times [0,T]).
\end{equation}  
Then we consider the second term in the right side of
(\ref{aux2}). Evidently,
$$
\Big\|\int_{0}^{t}
e^{(t-s)\{\sqrt{2\pi}\widehat{J}(p)+ibp+a\}}\sqrt{2\pi}p^{2}\widehat{G}(p)\widehat{f}_{v}(p,s)ds
\Big\|_{L^{2}({\mathbb R})}\leq 
$$
$$
\int_{0}^{t}\Big\|
e^{(t-s)\{\sqrt{2\pi}\widehat{J}(p)+ibp+a\}}\sqrt{2\pi}p^{2}\widehat{G}(p)\widehat{f}_{v}(p,s)
\Big\|_{L^{2}({\mathbb R})}ds.
$$
Note that
$$
\Big\|e^{(t-s)\{\sqrt{2\pi}\widehat{J}(p)+ibp+a\}}\sqrt{2\pi}p^{2}\widehat{G}(p)\widehat{f}_{v}(p,s)
\Big\|_{L^{2}({\mathbb R})}^{2}=
$$
\begin{equation}
\label{etsp2gpl2}
\int_{-\infty}^{\infty}
e^{2(t-s)\sqrt{2\pi}Re[\widehat{J}(p)]}e^{2a(t-s)}2\pi|p^{2}\widehat{G}(p)|^{2}|\widehat{f}_{v}(p,s)|^{2}
dp.
\end{equation}  
The right side of (\ref{etsp2gpl2}) can be bounded above using (\ref{fub1})
as
$$
2\pi e^{2aT}\|p^{2}\widehat{G}(p)\|_{L^{\infty}({\mathbb R})}^{2}
\int_{-\infty}^{\infty}|\widehat{f}_{v}(p,s)|^{2}dp\leq
$$
$$
e^{2aT}\Big\|\frac{d^{2}G}{dx^{2}}\Big\|_{L^{1}({\mathbb R})}^{2}
\|F(v(x,s),x)\|_{L^{2}({\mathbb R})}^{2}.
$$
By virtue of inequality (\ref{lubs}) we derive the estimate
$$
\Big\|e^{(t-s)\{\sqrt{2\pi}\widehat{J}(p)+ibp+a\}}\sqrt{2\pi}p^{2}\widehat{G}(p)\widehat{f}_{v}(p,s)
\Big\|_{L^{2}({\mathbb R})}\leq
$$
$$
e^{aT}\Big\|\frac{d^{2}G}{dx^{2}}\Big\|_{L^{1}({\mathbb R})}\{k\|v(x,s)\|_
{L^{2}({\mathbb R})}+\|h(x)\|_{L^{2}({\mathbb R})}\},
$$
so that
$$
\Big\|\int_{0}^{t}e^{(t-s)\{\sqrt{2\pi}\widehat{J}(p)+ibp+a\}}\sqrt{2\pi}p^{2}\widehat{G}(p)
\widehat{f}_{v}(p,s)ds\Big\|_{L^{2}({\mathbb R})}\leq
$$
$$
ke^{aT}\Big\|\frac{d^{2}G}{dx^{2}}\Big\|_{L^{1}({\mathbb R})}\int_{0}^{T}\|v(x,s)\|_
{L^{2}({\mathbb R})}ds+
Te^{aT}\Big\|\frac{d^{2}G}{dx^{2}}\Big\|_{L^{1}({\mathbb R})}\|h(x)\|_{L^{2}({\mathbb R})}.
$$
Recall upper bound (\ref{sch}). This gives us
$$
\Big\|\int_{0}^{t}e^{(t-s)\{\sqrt{2\pi}\widehat{J}(p)+ibp+a\}}\sqrt{2\pi}p^{2}\widehat{G}(p)
\widehat{f}_{v}(p,s)ds\Big\|_{L^{2}({\mathbb R})}^{2}\leq
$$
$$
e^{2aT}\Big\|\frac{d^{2}G}{dx^{2}}\Big\|_{L^{1}({\mathbb R})}^{2}
\{k\|v(x,s)\|_{L^{2}({\mathbb R}\times [0,T])}\sqrt{T}+
\|h(x)\|_{L^{2}({\mathbb R})}T\}^{2}.
$$
Hence, we obtain 
$$
\Big\|\int_{0}^{t}e^{(t-s)\{\sqrt{2\pi}\widehat{J}(p)+ibp+a\}}\sqrt{2\pi}p^{2}\widehat{G}(p)
\widehat{f}_{v}(p,s)ds\Big\|_{L^{2}({\mathbb R}\times [0, T])}^{2}\leq
$$
$$
e^{2aT}\Big\|\frac{d^{2}G}{dx^{2}}\Big\|_{L^{1}({\mathbb R})}^{2}
\{k\|v(x,s)\|_{L^{2}({\mathbb R}\times [0,T])}+
\|h(x)\|_{L^{2}({\mathbb R})}\sqrt{T}\}^{2}T^{2}<\infty
$$
under the given conditions for
$v(x,s)\in W^{(1, 2), 2}({\mathbb R}\times [0, T])$. Thus,
\begin{equation}
\label{int0tp2Gpfv}  
\int_{0}^{t}e^{(t-s)\{\sqrt{2\pi}\widehat{J}(p)+ibp+a\}}\sqrt{2\pi}p^{2}\widehat{G}(p)
\widehat{f}_{v}(p,s)ds\in L^{2}({\mathbb R}\times [0, T]).
\end{equation}
By means of (\ref{aux2}) along with statements (\ref{p2u0hpl2}), (\ref{int0tp2Gpfv}), we arrive at
\begin{equation}
\label{up2t12}  
p^{2}\widehat{u}(p,t)\in L^{2}({\mathbb R}\times [0, T]).
\end{equation}  
Therefore,
\begin{equation}
\label{d2udx2l2}
\frac{\partial^{2}u}{\partial x^{2}}\in L^{2}({\mathbb R}\times [0, T]).
\end{equation}  
It follows from (\ref{aux}) that
\begin{equation}
\label{duhdt}
\frac{\partial \widehat{u}}{\partial t}=\{\sqrt{2\pi}\widehat{J}(p)+ibp+a\}\widehat{u}(p,t)+
\sqrt{2\pi}\widehat{G}(p)\widehat{f}_{v}(p,t).
\end{equation}
According to (\ref{upt12}),
\begin{equation}
\label{aupt12}  
a\widehat{u}(p,t)\in L^{2}({\mathbb R}\times [0,T]).
\end{equation}
The norm can be easily estimated as
$$
\|ibp\widehat{u}(p,t)\|_{L^{2}({\mathbb R}\times [0,T])}^{2}=b^{2}\int_{0}^{T}
\Big\{\int_{|p|\leq 1}p^{2}|\widehat{u}(p,t)|^{2}dp+
\int_{|p|>1}p^{2}|\widehat{u}(p,t)|^{2}dp\Big\}dt\leq
$$
$$
b^{2}\{\|\widehat{u}(p,t)\|_{L^{2}({\mathbb R}\times [0,T])}^{2}+
\|p^{2}\widehat{u}(p,t)\|_{L^{2}({\mathbb R}\times [0,T])}^{2}\}<\infty
$$
due to (\ref{upt12}) and (\ref{up2t12}). Hence,
\begin{equation}
\label{bupt12} 
ibp\widehat{u}(p,t)\in L^{2}({\mathbb R}\times [0,T]).
\end{equation}
Recall inequality (\ref{fub}). Obviously,
$$
\|\sqrt{2\pi}\widehat{J}(p)\widehat{u}(p,t)\|_{L^{2}({\mathbb R}\times [0,T])}^{2}=2\pi
\int_{0}^{T}\int_{-\infty}^{\infty}|\widehat{J}(p)\widehat{u}(p,t)|^{2}dpdt\leq
$$
$$
2\pi\|\widehat{J}(p)\|_{L^{\infty}({\mathbb R})}^{2}\int_{0}^{T}\int_{-\infty}^{\infty}|\widehat{u}(p,t)|^{2}dpdt\leq
\|J(x)\|_{L^{1}({\mathbb R})}^{2}\|\widehat{u}(p,t)\|_{L^{2}({\mathbb R}\times [0,T])}^{2}.
$$
Then
\begin{equation}
\label{jhuhpt}
\sqrt{2\pi}\widehat{J}(p)\widehat{u}(p,t)\in L^{2}({\mathbb R}\times [0,T])
\end{equation}
via Assumption 1.1 along with (\ref{upt12}). 
Let us combine statements (\ref{aupt12}), (\ref{bupt12}) and (\ref{jhuhpt}). Thus,
\begin{equation}
\label{uptalfab12}
(\sqrt{2\pi}\widehat{J}(p)+ibp+a)\widehat{u}(p,t)\in L^{2}({\mathbb R}\times [0,T]).
\end{equation}
We turn our attention to deriving the upper bound on the norm of the
remaining term in the right side of (\ref{duhdt}) by virtue of (\ref{fub}) and 
(\ref{lubs}). Evidently,
$$
\|\sqrt{2\pi}\widehat{G}(p)\widehat{f}_{v}(p,t)\|_
{L^{2}({\mathbb R}\times [0,T])}^{2}\leq 2\pi
\|\widehat{G}(p)\|_{L^{\infty}({\mathbb R})}^{2}\int_{0}^{T}
\|F(v(x,t), x)\|_{L^{2}({\mathbb R})}^{2}dt\leq
$$
$$
\|G(x)\|_{L^{1}({\mathbb R})}^{2}\int_{0}^{T}(k\|v(x,t)\|_{L^{2}({\mathbb R})}+
\|h(x)\|_{L^{2}({\mathbb R})})^{2}dt\leq
$$
$$
\|G(x)\|_{L^{1}({\mathbb R})}^{2}\{2k^{2}\|v(x,t)\|_{L^{2}({\mathbb R}\times [0, T])}^{2}+
2\|h(x)\|_{L^{2}({\mathbb R})}^{2}T\}<\infty
$$
under the stated assumptions with 
$v(x,t)\in W^{(1, 2), 2}({\mathbb R}\times [0, T])$. Hence,
\begin{equation}
\label{ghfvhpt}  
\sqrt{2\pi}\widehat{G}(p)\widehat{f}_{v}(p,t)\in L^{2}({\mathbb R}\times [0,T]).
\end{equation}
By virtue of (\ref{duhdt}) along with statements (\ref{uptalfab12}) and
(\ref{ghfvhpt}),
$$
\frac{\partial \widehat{u}}{\partial t}\in L^{2}({\mathbb R}\times [0,T]).
$$
This means that
\begin{equation}
\label{dudtl2}  
\frac{\partial u}{\partial t}\in L^{2}({\mathbb R}\times [0,T])
\end{equation}
as well. Recall the definition of the norm (\ref{122n}).
According to (\ref{uxtl2}), (\ref{d2udx2l2}) and (\ref{dudtl2}), we have that
for the function uniquely determined by (\ref{aux}), 
$$
u(x,t)\in W^{(1, 2), 2}({\mathbb R}\times [0, T]).
$$
Thus, under the given conditions problem (\ref{aux}) defines a map
$$
\tau_{a, b}: W^{(1, 2), 2}({\mathbb R}\times [0, T])\to
W^{(1, 2), 2}({\mathbb R}\times [0, T]).
$$
The goal is to establish that under the our assumptions such map is a
strict contraction. Let us choose arbitrarily
$v_{1, 2}(x,t)\in W^{(1, 2), 2}({\mathbb R}\times [0, T])$. According to the
reasoning above,
$u_{1, 2}:=\tau_{a, b}v_{1, 2}\in W^{(1, 2), 2}({\mathbb R}\times [0, T])$.
It follows from (\ref{aux}) that
$$
\widehat{u_{1}}(p,t)=
$$
\begin{equation}
\label{1aux}
e^{t\{\sqrt{2\pi}\widehat{J}(p)+ibp+a\}}\widehat{u_{0}}(p)+\int_{0}^{t}
e^{(t-s)\{\sqrt{2\pi}\widehat{J}(p)+ibp+a\}}\sqrt{2\pi}\widehat{G}(p)\widehat{f}_{v_{1}}(p,s)ds,
\end{equation}
$$
\widehat{u_{2}}(p,t)=
$$
\begin{equation}
\label{2aux}
e^{t\{\sqrt{2\pi}\widehat{J}(p)+ibp+a\}}\widehat{u_{0}}(p)+\int_{0}^{t}
e^{(t-s)\{\sqrt{2\pi}\widehat{J}(p)+ibp+a\}}\sqrt{2\pi}\widehat{G}(p)\widehat{f}_{v_{2}}(p,s)ds.
\end{equation}
Here $\widehat{f}_{v_{j}}(p,s)$ with $j=1, 2$ denotes the Fourier image of
$F(v_{j}(x,s), x)$ under transform (\ref{ft}). From system
(\ref{1aux}), (\ref{2aux}) we easily deduce that
$$
\widehat{u_{1}}(p,t)-\widehat{u_{2}}(p,t)=
$$
\begin{equation}
\label{u1u2hint0t}  
\int_{0}^{t}
e^{(t-s)\{\sqrt{2\pi}\widehat{J}(p)+ibp+a\}}\sqrt{2\pi}\widehat{G}(p)
[\widehat{f}_{v_{1}}(p,s)-\widehat{f}_{v_{2}}(p,s)]ds.
\end{equation}
Clearly, the estimate on the norm 
$$
\|\widehat{u_{1}}(p,t)-\widehat{u_{2}}(p,t)\|_{L^{2}({\mathbb R})}\leq
$$
\begin{equation}
\label{int0tetsf12}  
\int_{0}^{t}\|
e^{(t-s)\{\sqrt{2\pi}\widehat{J}(p)+ibp+a\}}\sqrt{2\pi}\widehat{G}(p)
[\widehat{f}_{v_{1}}(p,s)-\widehat{f}_{v_{2}}(p,s)]\|_{L^{2}({\mathbb R})}ds
\end{equation}
holds.
Let us use (\ref{re}) along with inequality (\ref{fub}) to derive the upper bound as
$$
\|e^{(t-s)\{\sqrt{2\pi}\widehat{J}(p)+ibp+a\}}\sqrt{2\pi}\widehat{G}(p)
[\widehat{f}_{v_{1}}(p,s)-\widehat{f}_{v_{2}}(p,s)]\|_{L^{2}({\mathbb R})}^{2}=
$$
$$
2\pi \int_{-\infty}^{\infty}e^{2(t-s)\sqrt{2\pi}Re[\widehat{J}(p)]}e^{2(t-s)a}|\widehat{G}(p)|^{2}
|\widehat{f}_{v_{1}}(p,s)-\widehat{f}_{v_{2}}(p,s)|^{2}dp\leq 
$$
$$
2\pi e^{2aT}\|\widehat{G}(p)\|_{L^{\infty}({\mathbb R})}^{2}\int_{-\infty}^{\infty}
|\widehat{f}_{v_{1}}(p,s)-\widehat{f}_{v_{2}}(p,s)|^{2}dp\leq 
$$
$$
e^{2aT}\|G(x)\|_{L^{1}({\mathbb R})}^{2}\|F(v_{1}(x,s), x)-F(v_{2}(x,s), x)\|_
{L^{2}({\mathbb R})}^{2}.
$$
Recall inequality (\ref{lip}). Thus,
\begin{equation}
\label{lipl2}  
\|F(v_{1}(x,s), x)-F(v_{2}(x,s), x)\|_{L^{2}({\mathbb R})}\leq l
\|v_{1}(x,s)-v_{2}(x,s)\|_{L^{2}({\mathbb R})},
\end{equation}
so that
$$
\|e^{(t-s)\{\sqrt{2\pi}\widehat{J}(p)+ibp+a\}}\sqrt{2\pi}\widehat{G}(p)
[\widehat{f}_{v_{1}}(p,s)-\widehat{f}_{v_{2}}(p,s)]\|_{L^{2}({\mathbb R})}\leq
$$
\begin{equation}
\label{eatgvl2}  
e^{aT}l\|G(x)\|_{L^{1}({\mathbb R})}\|v_{1}(x,s)-v_{2}(x,s)\|_{L^{2}({\mathbb R})}.   
\end{equation}
By means of (\ref{int0tetsf12}) and (\ref{eatgvl2}), we arrive at
$$
\|\widehat{u_{1}}(p,t)-\widehat{u_{2}}(p,t)\|_{L^{2}({\mathbb R})}\leq
$$
$$
e^{aT}l\|G(x)\|_{L^{1}({\mathbb R})}\int_{0}^{T}
\|v_{1}(x,s)-v_{2}(x,s)\|_{L^{2}({\mathbb R})}ds.
$$
According to the Schwarz inequality, 
\begin{equation}
\label{sch2}  
\int_{0}^{T}\|v_{1}(x,s)-v_{2}(x,s)\|_{L^{2}({\mathbb R})}ds\leq
\sqrt{\int_{0}^{T}\|v_{1}(x,s)-v_{2}(x,s)\|_{L^{2}({\mathbb R})}^{2}ds}\sqrt{T},
\end{equation}
such that 
$$  
\|\widehat{u_{1}}(p,t)-\widehat{u_{2}}(p,t)\|_{L^{2}({\mathbb R})}\leq
$$
\begin{equation}
\label{u1hu2hl2}  
e^{aT}l\sqrt{T}\|G(x)\|_{L^{1}({\mathbb R})}\|v_{1}(x,t)-v_{2}(x,t)\|_
{L^{2}({\mathbb R}\times [0, T])}.
\end{equation}
This means that
$$
\|u_{1}(x,t)-u_{2}(x,t)\|_{L^{2}({\mathbb R}\times [0, T])}^{2}=\int_{0}^{T}
\|\widehat{u_{1}}(p,t)-\widehat{u_{2}}(p,t)\|_{L^{2}({\mathbb R})}^{2}dt\leq
$$
\begin{equation}
\label{u1u2l2v1v2}
e^{2aT}l^{2}T^{2}\|G(x)\|_{L^{1}({\mathbb R})}^{2}
\|v_{1}(x,t)-v_{2}(x,t)\|_{L^{2}({\mathbb R}\times [0, T])}^{2}.
\end{equation}  
Using (\ref{u1u2hint0t}), we obtain
$$
p^{2}[\widehat{u_{1}}(p,t)-\widehat{u_{2}}(p,t)]=
\int_{0}^{t}
e^{(t-s)\{\sqrt{2\pi}\widehat{J}(p)+ibp+a\}}\sqrt{2\pi}p^{2}\widehat{G}(p)
[\widehat{f}_{v_{1}}(p,s)-\widehat{f}_{v_{2}}(p,s)]ds.
$$
Obviously, the estimate from above on the norm 
$$
\|p^{2}[\widehat{u_{1}}(p,t)-\widehat{u_{2}}(p,t)]\|_{L^{2}({\mathbb R})}\leq
$$
\begin{equation}
\label{int0tetsf122}  
\int_{0}^{t}\|
e^{(t-s)\{\sqrt{2\pi}\widehat{J}(p)+ibp+a\}}\sqrt{2\pi}p^{2}\widehat{G}(p)
[\widehat{f}_{v_{1}}(p,s)-\widehat{f}_{v_{2}}(p,s)]\|_{L^{2}({\mathbb R})}ds
\end{equation}
is valid.
Let us recall condition (\ref{re}) along with inequality (\ref{fub1}). Hence,
$$
\|e^{(t-s)\{\sqrt{2\pi}\widehat{J}(p)+ibp+a\}}\sqrt{2\pi}p^{2}\widehat{G}(p)
[\widehat{f}_{v_{1}}(p,s)-\widehat{f}_{v_{2}}(p,s)]\|_{L^{2}({\mathbb R})}^{2}=
$$
$$
2\pi \int_{-\infty}^{\infty}e^{2(t-s)\sqrt{2\pi}Re[\widehat{J}(p)]}e^{2(t-s)a}|p^{2}\widehat{G}(p)|^{2}
|\widehat{f}_{v_{1}}(p,s)-\widehat{f}_{v_{2}}(p,s)|^{2}dp\leq 
$$
$$
2\pi e^{2aT}\|p^{2}\widehat{G}(p)\|_{L^{\infty}({\mathbb R})}^{2}\int_{-\infty}^{\infty}
|\widehat{f}_{v_{1}}(p,s)-\widehat{f}_{v_{2}}(p,s)|^{2}dp\leq 
$$
$$
e^{2aT}\Big\|\frac{d^{2}G}{dx^{2}}\Big\|_{L^{1}({\mathbb R})}^{2}
\|F(v_{1}(x,s), x)-F(v_{2}(x,s), x)\|_{L^{2}({\mathbb R})}^{2}.
$$
By virtue of formula (\ref{lipl2}), we derive
$$
\|e^{(t-s)\{\sqrt{2\pi}\widehat{J}(p)+ibp+a\}}\sqrt{2\pi}p^{2}\widehat{G}(p)
[\widehat{f}_{v_{1}}(p,s)-\widehat{f}_{v_{2}}(p,s)]\|_{L^{2}({\mathbb R})}\leq 
$$
\begin{equation}
\label{etsp2gpeat}
e^{aT}l\Big\|\frac{d^{2}G}{dx^{2}}\Big\|_{L^{1}({\mathbb R})}\|v_{1}(x,s)-v_{2}(x,s)\|_
{L^{2}({\mathbb R})}.
\end{equation}
Let us use (\ref{int0tetsf122}) along with (\ref{etsp2gpeat}) and (\ref{sch2}). This yields
$$
\|p^{2}[\widehat{u_{1}}(p,t)-\widehat{u_{2}}(p,t)]\|_{L^{2}({\mathbb R})}\leq
$$
\begin{equation}
\label{p2u1hu2hpt}  
e^{aT}\sqrt{T}l\Big\|\frac{d^{2}G}{dx^{2}}\Big\|_{L^{1}({\mathbb R})}
\|v_{1}(x,t)-v_{2}(x,t)\|_{L^{2}({\mathbb R}\times [0, T])},
\end{equation}
so that
$$
\Big\|\frac{\partial^{2}}{\partial x^{2}}[u_{1}(x,t)-u_{2}(x,t)]\Big\|_
{L^{2}({\mathbb R}\times [0, T])}^{2}=\int_{0}^{T}
\|p^{2}[\widehat{u_{1}}(p,t)-\widehat{u_{2}}(p,t)]\|_{L^{2}({\mathbb R})}^{2}dt\leq     $$
\begin{equation}
\label{u1u2l2v1v22}
e^{2aT}l^{2}T^{2}\Big\|\frac{d^{2}G}{dx^{2}}\Big\|_{L^{1}({\mathbb R})}^{2}
\|v_{1}(x,t)-v_{2}(x,t)\|_{L^{2}({\mathbb R}\times [0, T])}^{2}.
\end{equation}  
It easily follows from (\ref{u1u2hint0t}) that
$$
\frac{\partial}{\partial t}[\widehat{u_{1}}(p,t)-\widehat{u_{2}}(p,t)]=
$$
$$
\{\sqrt{2\pi}\widehat{J}(p)+ibp+a\}[\widehat{u_{1}}(p,t)-\widehat{u_{2}}(p,t)]+
\sqrt{2\pi}\widehat{G}(p)[\widehat{f}_{v_{1}}(p,t)-\widehat{f}_{v_{2}}(p,t)].
$$
Therefore,
$$
\Big\|\frac{\partial}{\partial t}[\widehat{u_{1}}(p,t)-\widehat{u_{2}}(p,t)]
\Big\|_{L^{2}({\mathbb R})}\leq a\|\widehat{u_{1}}(p,t)-\widehat{u_{2}}(p,t)\|_
{L^{2}({\mathbb R})}+
$$
$$      
|b|\|p[\widehat{u_{1}}(p,t)-\widehat{u_{2}}(p,t)]\|_
{L^{2}({\mathbb R})}+
\sqrt{2\pi}\|\widehat{J}(p)[\widehat{u_{1}}(p,t)-\widehat{u_{2}}(p,t)]\|_{L^{2}({\mathbb R})}+    
$$
\begin{equation}
\label{ddtu1hptu2hpt}
\sqrt{2\pi}\|\widehat{G}(p)[\widehat{f}_{v_{1}}(p,t)-\widehat{f}_{v_{2}}(p,t)]\|_
{L^{2}({\mathbb R})}.
\end{equation}  
Recall inequality (\ref{u1hu2hl2}). Thus, the first term in the right side of
(\ref{ddtu1hptu2hpt}) can be bounded from above by
\begin{equation}
\label{au1hu2hl2}  
a\nu e^{aT}\sqrt{T}l\|v_{1}(x,t)-v_{2}(x,t)\|_{L^{2}({\mathbb R}\times [0, T])}
\end{equation}
with $\nu$ introduced in (\ref{q}).
We estimate the norm as
$$
\|p[\widehat{u_{1}}(p,t)-\widehat{u_{2}}(p,t)]\|_{L^{2}({\mathbb R})}^{2}=
$$
$$
\int_{|p|\leq 1}p^{2}|\widehat{u_{1}}(p,t)-\widehat{u_{2}}(p,t)|^{2}dp+
\int_{|p|>1}p^{2}|\widehat{u_{1}}(p,t)-\widehat{u_{2}}(p,t)|^{2}dp\leq 
$$
$$  
\|\widehat{u_{1}}(p,t)-\widehat{u_{2}}(p,t)\|_{L^{2}({\mathbb R})}^{2}+
\|p^{2}[\widehat{u_{1}}(p,t)-\widehat{u_{2}}(p,t)]\|_{L^{2}({\mathbb R})}^{2}.
$$
Let us use (\ref{u1hu2hl2}) and  (\ref{p2u1hu2hpt}). This gives us
the upper bound on the second term in the right side of
(\ref{ddtu1hptu2hpt}) equal to
\begin{equation}
\label{bu1hu2hl2}  
|b|\nu e^{aT}\sqrt{T}l\|v_{1}(x,t)-v_{2}(x,t)\|_{L^{2}({\mathbb R}\times [0, T])}.
\end{equation}
By means of (\ref{fub}) and (\ref{u1hu2hl2}),
$$
2\pi\|\widehat{J}(p)[\widehat{u_{1}}(p,t)-\widehat{u_{2}}(p,t)]\|_{L^{2}({\mathbb R})}^{2}\leq 
2\pi \|\widehat{J}(p)\|_{L^{\infty}({\mathbb R})}^{2}\int_{-\infty}^{\infty}|\widehat{u_{1}}(p,t)-\widehat{u_{2}}(p,t)|^{2}dp\leq
$$
$$
\|J(x)\|_{L^{1}({\mathbb R})}^{2}\|\widehat{u_{1}}(p,t)-\widehat{u_{2}}(p,t)\|_{L^{2}({\mathbb R})}^{2}\leq
$$
$$
\|J(x)\|_{L^{1}({\mathbb R})}^{2}e^{2aT}l^{2}T\|G(x)\|_{L^{1}({\mathbb R})}^{2}\|v_{1}(x,t)-v_{2}(x,t)\|_{L^{2}({\mathbb R}\times [0, T])}^{2}.
$$
Therefore, the third term in the
right side of (\ref{ddtu1hptu2hpt}) can be estimated from above by
\begin{equation}
\label{bu1hu2hl3}  
\|J(x)\|_{L^{1}({\mathbb R})}\nu e^{aT}\sqrt{T}l\|v_{1}(x,t)-v_{2}(x,t)\|_{L^{2}({\mathbb R}\times [0, T])}.
\end{equation}
By virtue of  (\ref{fub}) along with inequality (\ref{lipl2}), we derive
$$
2\pi\int_{-\infty}^{\infty}|\widehat{G}(p)|^{2}
|\widehat{f}_{v_{1}}(p,t)-\widehat{f}_{v_{2}}(p,t)|^{2}dp\leq
$$
$$
2\pi
\|\widehat{G}(p)\|_{L^{\infty}({\mathbb R})}^{2}\int_{-\infty}^{\infty}
|\widehat{f}_{v_{1}}(p,t)-\widehat{f}_{v_{2}}(p,t)|^{2}dp\leq
$$
$$
\|G(x)\|_{L^{1}({\mathbb R})}^{2}\|F(v_{1}(x,t), x)-F(v_{2}(x,t), x)\|_
{L^{2}({\mathbb R})}^{2}\leq
$$
$$
\|G(x)\|_{L^{1}({\mathbb R})}^{2}l^{2}\|v_{1}(x,t)-v_{2}(x,t)\|_{L^{2}({\mathbb R})}^{2}.  
$$
This means that the fourth term in the right side of (\ref{ddtu1hptu2hpt}) can be bounded
from above by
\begin{equation}
\label{bu1hu2hl4}  
\nu l\|v_{1}(x,t)-v_{2}(x,t)\|_{L^{2}({\mathbb R})}.
\end{equation}
Let us combine (\ref{au1hu2hl2}), (\ref{bu1hu2hl2}), (\ref{bu1hu2hl3}) and
(\ref{bu1hu2hl4}). We obtain
$$
\Big\|\frac{\partial}{\partial t}[\widehat{u_{1}}(p,t)-\widehat{u_{2}}(p,t)]
\Big\|_{L^{2}({\mathbb R})}\leq
$$
$$
\nu e^{aT}\sqrt{T}l\{a+|b|+\|J(x)\|_{L^{1}({\mathbb R})}\}\|v_{1}(x,t)-v_{2}(x,t)\|_
{L^{2}({\mathbb R}\times [0, T])}+
$$
$$
\nu l\|v_{1}(x,t)-v_{2}(x,t)\|_{L^{2}({\mathbb R})}.
$$
This allows us to estimate the norm as
$$
\Big\|\frac{\partial}{\partial t}(u_{1}(x,t)-u_{2}(x,t))\Big\|_
{L^{2}({\mathbb R}\times [0, T])}^{2}=\int_{0}^{T}
\Big\|\frac{\partial}{\partial t}[\widehat{u_{1}}(p,t)-\widehat{u_{2}}(p,t)]
\Big\|_{L^{2}({\mathbb R})}^{2}dt\leq           
$$
\begin{equation}
\label{ddtu1u2l2}  
2\nu^{2}l^{2}[e^{2aT}T^{2}\{a+|b|+\|J(x)\|_{L^{1}({\mathbb R})}\}^{2}+1]
\|v_{1}(x,t)-v_{2}(x,t)\|_{L^{2}({\mathbb R}\times [0, T])}^{2}.
\end{equation}
Recall the definition of the norm (\ref{122n}). Using upper
bounds (\ref{u1u2l2v1v2}), (\ref{u1u2l2v1v22}) and (\ref{ddtu1u2l2}),
we arrive at
$$
\|u_{1}-u_{2}\|_{W^{(1, 2), 2}({\mathbb R}\times [0, T])}\leq
$$
\begin{equation}
\label{contr}
\nu l\sqrt{T^{2}e^{2aT}(1+2[a+|b|+\|J(x)\| _{L^{1}({\mathbb R})}]^{2})+2}
\|v_{1}-v_{2}\|_{W^{(1, 2), 2}({\mathbb R}\times [0, T])}.
\end{equation}
The constant in the right side of (\ref{contr}) is less than one
via (\ref{qlt}). Therefore, under the stated assumptions equation
(\ref{aux}) defines the map
$$
\tau_{a, b}: W^{(1, 2), 2}({\mathbb R}\times [0, T])\to W^{(1, 2), 2}
({\mathbb R}\times [0, T]),
$$
which is a strict contraction. Its unique fixed
point $w(x,t)$ is the only solution of
problem (\ref{h}), (\ref{ic}) in
$W^{(1, 2), 2}({\mathbb R}\times [0, T])$. \hfill\lanbox

\bigskip

\noindent
{\it Proof of Corollary 1.4.} The validity of the statement of the Corollary
comes from the fact that the constant in the right side of inequality
(\ref{contr}) is independent of the initial condition (\ref{ic}) (see e.g. ~\cite{EZ01}). Therefore,
problem (\ref{h}), (\ref{ic}) admits a unique solution
$w(x,t)\in W^{(1, 2), 2}({\mathbb R}\times {\mathbb R}^{+})$. Suppose
that $w(x,t)$ is trivial  for all $x\in {\mathbb R}$ and $t\in {\mathbb R}^{+}$.
This will give us the contradiction to the condition that
$\hbox{supp}\widehat{F(0, x)}\cap \hbox{supp}\widehat{G}$ is a set of
nonzero Lebesgue measure on the real line. \hfill\lanbox

\bigskip

%%%%%%%%%%%%%%%%%%%%%%%%%%%%%%%%%%%%%%%%%%%%%%%%%%%%%%%%%%%%%%%%

\centerline{\bf 3. Acknowledgement}

\bigskip

\noindent
V.V. is grateful to Israel Michael Sigal for the partial support by the
NSERC grant NA 7901.

\end{document}